\documentclass[10pt,a4paper]{article}
\usepackage{amsmath,amssymb,enumerate,verbatim}
\usepackage{amsthm}
\usepackage{epsfig,fancyhdr,color}
\usepackage{hyperref}
\usepackage{makeidx}

\usepackage[utf8]{inputenc}
\usepackage[T1]{fontenc}

\makeindex

\newcommand{\ds}{\displaystyle}
\newcommand{\tq}{\, \big| \, }

\DeclareMathOperator{\met}{Met}
\DeclareMathOperator{\aire}{Area}
 \renewcommand{\r}{\mathbb{R}}

\newtheorem{theorem}{\rm\bf Theorem}[section]
\newtheorem{proposition}[theorem]{\rm\bf Proposition}
\newtheorem{lemme}[theorem]{\rm\bf Lemma}
\newtheorem{corollaire}[theorem]{\rm\bf Corollary}

\theoremstyle{definition}
\newtheorem{definition}[theorem]{\rm\bf Definition}

\theoremstyle{remark}

\date{\today}

\title{On Alexandrov's Surfaces with Bounded Integral Curvature}

\author{Marc Troyanov, EPFL\footnote{EPFL, Institut de Math\'{e}matiques,    1015 Lausanne, Switzerland,   marc.troyanov@epfl.ch}}

\begin{document}

\maketitle

\maketitle

\begin{abstract}
During  the years 1940-1970, Alexandrov and the ``Leningrad School''  have  investigated  the geometry of singular surfaces in depth. 
The theory developed by this school is about  topological surfaces with an intrinsic metric for  which we can define a notion of curvature, which is a Radon measure.  This class of surfaces has good convergence properties   and   is remarkably stable with respect  to various  geometrical constructions (gluing etc.). It includes  
polyhedral surfaces as well as Riemannian surfaces of class   $C^2$, and both of these classes are dense families of  Alexandrov's surfaces. 
Any singular surface that can be  reasonably thought of is an Alexandrov surface and  
a number of  geometric properties of smooth surfaces extend and  generalize  to this class. 

The goal of this paper is to give an introduction to   Alexandrov's theory, to provide some examples and state some of the  
fundamental facts of the theory. We discuss the conformal viewpoint introduced by Yuri G. Reshetnyak and explain how it leads to a 
classification of compact  Alexandrov's surfaces. 

\medskip

This article is an updated and translated version of  the paper \cite{Troyanov2009}, to be included in the forthcoming  book
\emph{Reshetnyak's Theory of Subharmonic Metrics}  edited by Fran\c{c}ois Fillastre and Dmitriy Slutskiy and to
be published by Springer and the  Centre de recherches mathématiques (CRM) in Montr\'eal. 

\medskip

\noindent AMS Mathematics Subject Classification:     \  53.c45, 52.b70
\end{abstract}


\bigskip


\section{Introduction}  
The mathematical object that we call a ``surface''  is not unequivocally defined, and the choice of one particular mathematical framework rather than another one may lead us 
to study quite different objects such as smooth surfaces,  polyhedral surfaces, (semi-)algebraic surfaces etc.   A famous story about Henri Lebesgue says that, following a class by Gaston  Darboux, where the master proved that any developable surface  in the Euclidean three-space (that is a surface locally isometric to the plane)  is a ruled surface
(that is through every point on the surface there is a straight segment contained in the surface),  the young Lebesgue 
threw a  crumpled handkerchief on the table and said  \emph{``this handkerchief is a developable surface, yet it is not ruled''}. Darboux was doing classical differential geometry and Lebesgue observed that, in nature, surfaces are rarely smooth.

\medskip

Can we \emph{unify}, rather than \emph{oppose}, Darboux's smooth differential geometry and Lebesgue's crumpled geometry? Is there a non-trivial theory that encompasses all the surfaces we can reasonably conceive of?  In this article, I will give a rough sketch of the  solution that Alexandre Danilovitch Alexandrov and his collaborators gave to this problem during the years 1940-1970.

\medskip

Before moving on, let us start with an easy to state Problem. We will call \emph{smooth surface} a 2-dimensional  Riemannian manifold $(S,g)$, and \emph{polyhedral surface}  a metric space which is locally isometric to a 2-dimensional polyhedron (convex or not). We then ask:

\newpage

\textbf{Problem 1. }
\begin{enumerate}[(a)]
  \item Is any smooth surface  the limit of a sequence of polyhedral surfaces?  
  \item Is any polyhedral surface the limit of a sequence of smooth surfaces?
  \item What kind of geometric invariants  pass to the limit for these convergences? 
\end{enumerate}

\smallskip 

The reader will have noticed that these questions are not well formulated: what kind of convergence are we talking about? But let us ignore this point  for the time being. 

Problem (1c) is of fundamental importance. If our aim is to be able to move from one type of surface to another, from the smooth world to the polyhedral world, it is necessary 
to understand which geometric invariants are stable under convergence. This is a hard problem, and the first empirical observations are not too inviting. 

\smallskip 

\begin{enumerate}[$\circ$]
\item[$\circ$] The first example is the \emph{Schwarz Lantern} \cite{Schwarz}, which is a well-known example of a sequence of (non-convex) polyhedra $P_i$  converging to a cylinder in Euclidean three-space. The area $A$ of this cylinder is finite but the area of $P_i$ tends towards $+\infty$ when  $i\to \infty$ (more generally we can adjust the parameters of the polyhedra to obtain a  sequence $\{P_i\}$  whose area converges to any preassigned  number in the range $[A,+\infty]$, see \cite [p. 117]{Gugg}). This, and other related phenomena, are studied in the nice paper  \cite{HPW}.

\item[$\circ$] A second example is given by  the group of isometries. Consider a sequence of polyhedra $P_i$ in Euclidean three-space converging to the round  sphere. The group $G_i$ of isometries of $P_i$ is a finite subgroup of $O(3)$; let  us assume that the $P_i$'s are increasingly symmetric, that is $G_i \subset G_{i+1}$, and denote by $G$ the Hausdorff limit of the sequence $\{G_i\} \subset O(3)$ (that is the closure of $\cup_i G_i$). Then $G$ is either a finite group, or it contains a finite index subgroup conjugate to $SO(2)$. In particular it is a Lie group of dimension $0$ or $1$, yet the group of isometries of the limit sphere is the three dimensional Lie group $O(3)$.
\end{enumerate}
 
\smallskip

\textbf{Exercice.} Describe all subgroups of $O(3)$, up to conjugation, that are limits of a sequence
of finite subgroups (there are five types of such groups, named after Pierre Curie\footnote{There are in fact 7 types of Curie's groups, but that list includes $SO(3)$ and the full group $O(3)$. Only five types of Curie's groups are limits of finite groups.}, see \cite{Shubnikov}).

\bigskip

These examples should  convince us that investigating the geometric convergence of surfaces quickly leads us to 
some  deep  geometrical questions. It will be useful to formalize  our Problem 1 in a more precise way:

\medskip

\textbf{Problem 1'}. Given a topological surface $S$, we wish  to define a space $\mathcal{M}(S)$ containing all ``reasonable''  metrics on $S$ and provide this space with a topology for which both the polyhedral metrics and the Riemannian metrics on $S$ form dense subsets. 
We also wish  to describe those geometric invariants that define continuous functions on $\mathcal{M}(S)$.

\medskip

Surfaces with bounded integral curvature in the sense of Alexandrov prove to be an adequate answer to this problem, and it is our goal in this paper to give 
a brief introduction to the subject, mainly from the conformal viewpoint developed by Y. Reshetnyak. 

\medskip

Note that a number of results in the  Alexandrov's theory of surfaces with bounded integral curvature have been published in Russian and have not  been translated. The books \cite{AZ} and \cite{R93} as well as as the current volume are useful references in English. The subject is currently quite active, we refer to  \cite{Adamowicz,Barthelme,BB1,BB2,BL,CreutzRomney,Debin,Richard2018,Shioya} for some 
recent works on the subject. One may also expect some fruitful relations with  the intense developments of discrete differential geometry such as exposed in the books  \cite{Bobenko} and \cite{GuYau}.


\section{Definition of Alexandrov's surfaces}

We now define the central notion of this paper: Let $S$ be a  fixed  topological surface.  For simplicity  we assume, throughout the paper, the surface $S$  to be closed and oriented.

\medskip

\begin{definition}
A \textit{metric with bounded integral curvature in Alexandrov's sense}  \index{surfaces with bounded integral curvature}
on $S$ is  given by a continuous function $d: S \times S \to \mathbb{R}$ such that 
\begin{enumerate}[(i)]
  \item $d$ is a metric  that induces the manifold topology on $S$.
  \item The distance $d$ is \textit{intrinsic}, that is for any pair of points $x,y \in S$, there is a sequence of continuous curves $\gamma_n : [0,1] \to S$ joining $x$ to $y$ such that the $d$-length of $\gamma_n$ converges to $d(x,y)$. 
\item The distance $d$ on $S$ is the uniform limit of the distances  $d_i$ associated to a sequence of
Riemannian metrics $g_i$ on $S$ such that the integral of the absolute value of the curvature is uniformly bounded.
\end{enumerate}
\end{definition}

\medskip

This definition calls for some explanations. 
Recall first that the \textit{length} of  a curve $\gamma : [a,b] \to S$ with respect to a metric $d$ is defined as
$$
 \ell (\gamma)  = \sup \sum_{k=0}^m d(\gamma (t_k),\gamma (t_{k+1})),
$$
where the supremum is taken over all subdivisions $a = t_0, < t_1 < \cdots < t_m = b$. Condition (ii) in the above definition then says that $(S,d)$ is
a \emph{length space}, that is the distance between any pair of points is the infimum of the lengths of all curves joining them.

\medskip

The third condition justifies the name ``bounded integral curvature''. To clarify its meaning, we need some reminders and preliminaries.

\smallskip

Let us denote by $\met(S)$ the set of metrics satisfying the conditions (i) and (ii) in the above definition and $U \subset S$
an arbitrary subset of $U$. 
We then define the \textit{uniform distance} in  $U$ between two elements $d_1,d_2 \in \met(S)$  as
\begin{equation}
  D_U(d_1,d_2) = \sup \{ |d_1(x,y) - d_2(x,y) | : \, x,y \in U \}.
\end{equation}
The metric $d\in \met(S)$ is said to be \emph{Riemannian} if a smooth structure on $S$ and a Riemannian metric $g$ on $S$ are given, such that 
$$
 d(x,y) = \inf \int_0^1 \sqrt{g(\dot \gamma (t),\dot \gamma (t))} dt ,
$$
where the infimum is taken over all piecewise $C^1$ curves $\gamma: [0,1] \to S$ satisfying $\gamma (0) = x$ and $\gamma (1) = y$.

\medskip

Let us now recall the notions of area and curvature of $g$. On any simply connected domain $U$ of $S$, we can find  two  differential  forms $\theta^1, \theta^2 \in \Omega^1(U)$ of degree $1$ such that $\theta^1 \wedge \theta^2$ is
positive and 
$$
 g = (\theta^1)^2 + (\theta^2)^2 = \theta^1\otimes \theta^1 + \theta^2\otimes \theta^2.
$$
Such a pair of $1$-forms  $\theta^1, \theta^2$  is called a (positively oriented) \emph{orthonormal moving coframe}.

\medskip

The \emph{Hodge star}, is the unique linear map  $*: \Omega^k(S) \to \Omega^{2-k}(S)$, where $\Omega^k(S)$ is the space of differential forms of degree $k$ on $S$ ($k=0,1,2$), such that  the conditions 
$$
* (1) = \theta^1 \wedge \theta^2, \qquad * (\theta^1 \wedge \theta^2) = 1, \qquad * \theta^1 =\theta^2, \qquad * \theta^2 = -\theta^1
$$
hold for any positively oriented orthonormal moving coframe $\theta^1,\theta^2$. 
The \emph{connection form} associated with the moving coframe $\theta^1, \theta^2$ is the $1$-form $\omega \in \Omega^1(U)$ defined by
\begin{equation}\label{eq.formconnection}
 \omega = -(*d\theta^1)\theta^1 -(*d\theta^2)\theta^2,
\end{equation}
and it  is the unique  $1$-form on $U$ satisfying \textit{Elie Cartan's Structure Equations}:
$$
  \begin{cases}
  d\theta^1 &= - \omega \wedge \theta^2 \\
  d\theta^2 &= \omega \wedge \theta^1.
  \end{cases}
$$
A direct calculation gives us the following 
\begin{lemme}
The differential forms $dA_g = \theta^1 \wedge \theta^2$ and $d\omega$ are independent on the chosen positively oriented orthonormal moving coframe  $\theta^1, \theta^2 \in \Omega^1(U)$. They are thus globally defined two-forms on $S$.
\end{lemme}

\medskip

We have just  defined two  differential forms of degree $2$ on $S$.  Although they are not exact forms, we will denote them by  $d\omega$ and 
$
 dA_g = \theta^1 \wedge \theta^2.
$
We think of these two-forms as measures on $S$ (the measure $d\omega$ is a signed measure and $dA_g$ is  the area measure, that is the two-dimensional Hausdorff measure on $(S,d)$).

\medskip

\begin{definition}
We  call $dA_g$ the \textit{area measure} and $d\omega$ the \textit{curvature measure} of  $(S,g)$.
\end{definition}

 Note that the  curvature  measure $d\omega$ is clearly invariant by homothety.  The function $K : S \to \mathbb{R}$ defined by $K = *d\omega$ (i.e. $d\omega = K dA_g$) is the \emph{Gaussian curvature} of $g$ (the Gaussian curvature is thus the Radon Nikodym derivative of $d\omega$ with respect to $dA_g$).  

\medskip

The curvature measure of a  Borel set $E \subset S$ will be denoted by
$\ds \omega(E) = \int_E d\omega = \int_E kdA,$ and the Gauss-Bonnet formula tells us that the integral of the curvature of a closed surface 
$(S,g)$ is equal to $2\pi$ times  its Euler characteristic $\chi (S)$:
\begin{equation}\label{eq.GB}
 \omega (S)  = 2\pi \chi (S).
\end{equation}
It is convenient to also introduce the following (non negative) measures:
$$
 d\omega^+ = K^+dA_g, \qquad d\omega^- = K^-dA_g \quad \text{and}  \quad  d |\omega| = |K|dA_g,
$$
where $K^+ = \max \{K,0\}$ and  $K^- = \max \{-K,0\}$. Notice  that 
$d\omega = d\omega^+ -  d\omega^-$ and $d |\omega| =d\omega^+ + d\omega^-$.

\medskip

Having  these notions in mind, we can now complete our definition of surfaces with bounded integral curvature.
Condition (iii) says that there exists a sequence of Riemannian metrics $\{g_j\}$ on $S$ and a constant $C$  such that 
$$
 D_U(d_{g_j},d) \to 0  \quad  \mathrm{and}   \quad   |\omega_{g_j}|(S) =   \int_S |K_{g_i}|dA_i  \leq C, 
$$
for all $j\in  \mathbb{N}$.

\medskip 

From the Banach-Alaoglu Theorem  we know that there exists a Radon measure $d\omega$ on $S$ and a  subsequence $\{g_{j'}\}$   of $\{g_j\}$ such that  $\{d\omega_{g_{j'}}\}$  converge weakly to $d\omega$. 
Recall that $d\omega_{g_{j'}}$  converges  weakly\footnote{some authors prefer to say $*$-\textit{weakly}.} to $d\omega$ on $S$  if \ 
$$\ds \int_S f d\omega_{g_{j'}}  \to \int_S f d\omega$$ for any bounded continuous function $f\in C(S)$. This condition implies that $\omega_{g_{j'}}(E) \to \omega(E)$ for any measurable set $E\subset S$ whose frontier  $\mathrm{Fr}(E)$ satisfies $\omega(\mathrm{Fr}(E))=0$, see e.g.  \cite[p. 134]{Doob}.
A basic result is then the following 

\begin{theorem}\label{th.defomega}
If $(S,d)$ is an Alexandrov surface, then the  Radon measure $d\omega = \lim_{j'} d\omega_{g_{j'}}$ is well defined on  $S$. In particular it is independent of the  chosen sequence  of Riemannian metrics $\{g_j\}$.
\end{theorem}

This result is a consequence of Theorem \ref{th.conv.ALex} below. It 
 is a non trivial fact since only uniform convergence is assumed for the distance while the 
Riemannian curvature depends on the derivatives of the metric tensor $g$  up to order $2$.

\medskip

\begin{definition}
The  limit measure $d\omega $ is called the \emph{curvature measure}  \index{curvature measure} of the Alexandrov surface $(S,d)$.
\end{definition}

\medskip

\textbf{Remarks.}  \textbf{(1)}  The geometric meaning of   Theorem \ref{th.defomega} is that the curvature measure $d\omega$ of an Alexandrov surface  $(S,d)$ can be synthetically  constructed from  the metric $d$ alone, without referring to a sequence of  Riemannian approximations. The construction is given in   \cite[chap. 5]{AZ},
it is rather deep and technical but the basic geometric intuition  is to triangulate the surface with small geodesic triangles and to sum the 
\textit{angular excess} of each triangle\footnote{One should think of the angular excess of a triangle as the equal to (sum of the angles - $\pi$), 
but the notion of angle is  a priory not well defined and Alexandrov replaces it with a notion of  ``upperangle''.}, 
then let the mesh of the triangulation goes to $0$. 
 
\medskip

\textbf{(2)}  For simplicity, we have assumed the surface $S$ to be compact. In the non compact case, one only assumes condition (iii)
in the definition to hold locally. However one can no longer apply the Banach-Alaoglu Theorem and the existence of the limit 
measure $d\omega$ is more delicate to establish. The relevant measure theoretic development are exposed in  \cite{Alexandrov1943},
see also  in \cite[pages 233--250]{Alexandrov}.

\medskip

 \textbf{(3)}  The Gauss-Bonnet formula still holds for any  closed Alexandrov surface because of the continuity of the integrals with respect to   weak convergence of the measures. We will see  in Section \ref{sec.class} that conversely \emph{any Radon measure on a closed surface $S$ that is compatible with the Gauss-Bonnet formula an
such that $\omega(\{x\}) < 2\pi$ for any point $x\in S$  is the curvature   measure of  an Alexandrov metric on} $S$.

\medskip

 \textbf{(4)} Note that our definition of surfaces with bounded integral curvature in the sense of Alexandrov  is not the original definition given by Alexandrov in \cite{AZ}, and it is therefore a fundamental theorem of the theory. 
Alexandrov defines the surfaces with bounded integral curvature in a synthetic (that is purely metric) way.
The condition is  that the total angular excess of any family of non overlapping simple triangles in $S$ is  uniformly bounded. This definition is explained in details on pages 4--6 of \cite{AZ}, see also \cite [p. 71-74]{R93}.  The equivalence between our definition of compact surfaces with bounded integral curvature and the definition given in \cite{AZ} can be proved using the results of Chapter 3 in \cite{AZ}, where it is a proved than any compact Alexandrov surface is the uniform limit of a sequence of polyhedral surfaces with a uniform bound on the total absolute curvature,
together with  the (quite obvious) fact that any polyhedral surfaces is the uniform limit of a sequence of Riemannian surface.

\medskip

 \textbf{(5)}  Let us finally mention that Alexandrov proved in  \cite[page 240]{AZ} that the 
curvature measure of a family of Alexandrov Surfaces continuously depends on the metric. More precisely :

\begin{theorem} \label{th.conv.ALex}
Let  $(S,d)$ be a  compact surface with bounded integral measure in Alexandrov's sense and $\{ d_j\}$ a sequence of metrics with bounded integral curvature 
which uniformly converges to $d$. Suppose that $|\omega_{j}|(S) \leq C$ for some constant $C$, where $d|\omega_{j}|$ is the absolute curvature measure of 
the Alexandrov surface  $(S,d_j)$. 
Then the curvature measure of $(S,d)$ is the weak limit of the curvature measures of  $\{(S,d_j)\}$, that is 
$$
  D(d_j,d) \to 0 \quad \Rightarrow  \quad  \lim_{i\to \infty}\int_S f(x) d\omega_j (x) = \int_S f(x) d\omega (x)
$$
for any continuous function $f$ on $S$,  where $d\omega_j$ is the curvature measure of $(S,d_j)$ and 
 $d\omega$  is the curvature measure of $(S,d)$.
\end{theorem}

\section{On the Weyl Problem.} \label{WeylPb}

The work of Alexandrov was in part motivated by the following problem posed by H. Weyl in 1916: \textit{Is any smooth  Riemannian metric with positive curvature on the $2$-sphere isometric to the boundary of a unique convex domain in $\r^3$} ?
\index{Weyl's problem}

\medskip 

Note that the boundary of an arbitrary bounded convex domain in  $\r^3$ can be approximated by smooth convex surfaces and is thus a surface with bounded integral curvature. Alexandrov gave the following answer to the Weyl problem in 1948:

\begin{theorem}\label{th.Realisation.Alexandrov}
Every metric with non negative bounded integral curvature on the two dimensional sphere can be realized as the boundary of a bounded convex domain  $\Omega \subset  \r^3$.  
\end{theorem}

The proof can be found in  \cite[page 269]{Alexandrov}, Alexandrov also proved a version of this Theorem for convex domains in 
the three dimensional sphere or hyperbolic space, see chapter XII in the same reference. The proof  is a \textit{tour de force} based on the simple observation that polyhedral metrics on a closed surface $S$ (that is Alexandrov metrics on a surface with discrete curvature measure) are dense in the space of all Alexandrov metrics on $S$, combined with classical results on the rigidity of convex polyhedron and the convergence Theorem \ref{th.conv.ALex}. 

\medskip

Notice that Weyl formulated his problem for smooth surfaces. The regularity of the convex domain in the previous Theorem was later   proved   by A. V. Pogorelov, who also proved the uniqueness of 
the convex domain  up to an isometry of $\r^3$. An excellent reference on the Weyl Problem is the book \cite{BZ1993} by Yu. D. BuragoV. A. Zalgaller. We also mention the
recent papers  \cite{Prosanov} and  \cite{Schlenker} by Roman Prosanov and  Jean-Marc Schlenker  for an account of recent  developments of the subject.

\section{Alexandrov Surfaces obtained by Gluing Riemannian Pieces}

A  simple recipe  to construct  Alexandrov surfaces is by gluing  triangles. Consider a finite system $\{ T_1, T_2, \dots , T_m\}$  of Riemannian triangles, that is each $T_i$ is isometric  to a two dimensional disk  equipped with a  Riemannian metric such that the boundary $\partial T_i$ is piecewise $C^2$ with exactly  three angular points. 
Let us  glue these triangles according to a prescribed triangulation of a given surface  $S$, by identifying the matching side pairs by isometries. We then obtain a surface homeomorphic to the given surface  $S$ with a natural metric induced by the Riemannian metric in each triangle. 

This metric  space is an Alexandrov surface whose  curvature  measure is  given by 
$$
 d\omega = d\omega_0 + d\omega _1 + d\omega_2,
$$
where $d\omega_2$ is absolutely continuous with respect to the area measure  and given by
$$
 d\omega_2 = K dA
$$
inside each triangle. The measure $d\omega_0$ is a discrete measure supported by the vertices of the triangulation, and such that for each vertex $p$, we have
$$
 \omega_0(\{p\}) = 2\pi - (\text{the sum of the angles of all triangles $T_i$ incident with the vertex $p$}).
$$
The point $p$ is then a \textit{conical singularity}     \index{conical singularity}  if $\omega_0(\{p\})\neq0$ and  $\omega_0(\{p\}) < 2 \pi$ and a \textit{cusp}   \index{cusp}  (more precisely a \textit{finite cusp})  if $ \omega_0(\{p\}) = 2 \pi$. A cusp occurs at a point $p$ if and only $p$ is a vertex of the triangulation and all triangles incident with $p$ have a vanishing angle  at this point. 
 
\medskip

Finally the measure $d\omega _1$ is supported by the edges of the triangulation. On each edge $e = T_i\cap T_j$, we have
$$
 d\omega_1 = (k^+ - k^-) ds
$$
where $k^+$ and $k^-$ are the geodesic curvatures of the edge  $e$ seen in  both  triangles adjacent to $e$,
taking into account the edge orientation.

\medskip

A \emph{polyhedral surface} is a surface obtained by gluing   together flat Euclidean triangles (in the classical sense, that is triangle in the Euclidean plane whose edge are straight segments).   The curvature measure  is then concentrated at the vertices of the triangulation and in the neighborhood of each vertex, the surface is locally isometric to a Euclidean cone. For this reason a polyhedral surface is also called an \emph{Euclidean surface with conical singularities}. Conversely, any Alexandrov surface with a discrete curvature measure is a polyhedral surface.

\medskip

Let us look at some examples: Consider first  the surface of a \textit{cube} (which we can triangulate any way we want). The six faces of the cube are flat and therefore $d\omega_2 = 0$, the twelve edges  of the cube are geodesic as seen from each incident face and therefore $d\omega_1 = 0$. The eight vertices are each incident with  3 angles of $\frac{\pi}{2}$, therefore we have
$$
 \omega (p) = 2\pi - 3\frac{\pi}{2} = \frac{\pi}{2},
$$
at  each of the eight  vertices. The Gauss-Bonnet Formula for the cube is then the  simple fact  that 
$$
 \int_S d\omega = \int_S d\omega_0 = 8\frac{\pi}{2} = 4 \pi = 2\pi \chi(S^2).
$$
As a second example, consider a \textit{can}, that  is an Euclidean cylinder of radius $r$ and height $h$ together with  its top and bottom  which are  both  Euclidean discs $D_1,D_2$ of radius $r$. One can topologically triangulate the different parts of this can as desired. The cylinder and the two discs are flat surfaces, so the curvature is concentrated on the two circles bounding the two discs. These circles are geodesic as curves in the cylinder, and as curves in the discs $D_i$, they have   constant geodesic curvature $k = \frac{1}{r}$.
We therefore have $d\omega_0= d\omega_2 = 0$ and 
$$
  d\omega_1 = \frac{1}{r} \left. ds \right|_{\partial D_1} + \frac{1}{r} \left. ds \right|_{\partial D_2}.
$$
The Gauss-Bonnet Formula in this case says that 
$$
 \int_S d\omega = \int_S d\omega_1 = \frac{1}{r} \text{Length} (\partial D_1) +  \frac{1}{r} \text{Length} (\partial D_2)
 =  \frac{1}{r} (2\pi r + 2\pi r) = 4 \pi.
$$

Alexandrov Surfaces obtained by gluing Riemannian pieces are also discussed in  \cite{Strichartz}; this paper also discuss the spectrum of the Laplacian on such surfaces.

\section{Conformal Structure and Uniformisation of Smooth Surfaces}

In this section, we consider smooth Riemannian metrics and their conformal deformations. Recall that a Riemannian metric $\tilde{g}$ on the differentiable surface $S$ is a \emph{conformal deformation} of $g$ if there is a function: $u: S \to \r$ such that $\tilde{g} = e^{2u}g$. If $u$ is a constant, the metric $\tilde{g}$ is said to be \emph{homothetic} to $g$.  There is a simple formula relating the curvatures of two conformally  equivalent metrics. This formula involves the Laplacian operator on the Riemannian surface and we first recall its definition:

\medskip 

\begin{definition}
The \emph{Laplacian} (or \textit{Laplace-Beltrami operator}) of  $u$  with respect to the metric $g$  is the differential operator defined as
$$\Delta_g u = -* d*du.$$
One can equivalently write $\Delta_g u dA = -d*du$. 
\end{definition}
The Laplacian  is an elliptic operator, in local coordinates we have 
$$
 \Delta_g u = -\frac{1}{\sqrt{\det (g_{ij})} } \sum_{\mu,\nu = 1}^2 \frac{\partial}{\partial x_{\mu}}\left(g^{\mu\nu} {\sqrt{\det (g_{ij})} } \; \frac{\partial u}{\partial x_{\nu}}\right).
$$

\medskip

\begin{lemme}\label{lemconformalcurvature}
Suppose that   $g$ and  $\tilde{g} = e^{2u}g$ are two conformally equivalent smooth Riemannian metrics on a surface $S$, then 
\begin{equation}\label{eq.chgt.curvature1}
    d\widetilde{\omega} = d\omega  + \Delta_g u dA_g,
\end{equation}
where $d\widetilde{\omega}$ and $d\omega$ are  the curvature measures of $\tilde{g}$ and $g$, $dA_g$ is the area measure of $g$
and $\Delta_g$ is  its Laplacian.
\end{lemme} 

Equation \eqref{eq.chgt.curvature1} can also be written as
\begin{equation}\label{eq.chgt.curvature}
   \tilde{K}e^{2u} = K + \Delta_g u.
\end{equation}

\medskip 

\textbf{Proof.}  
Consider  an open domain $U \subset S$  in which a moving coframe  $\theta^1,\theta^2\in \Omega^1(U)$ such that $g = (\theta^1)^2 + (\theta^2)^2$ is given. Then $\tilde{\theta}^1 = e^u \theta^1$, $\tilde{\theta}^2 = e^u \theta^2$ is a moving coframe for the metric $\tilde{g} = e^{2u}g$ and using the structure equations, it is  not difficult to check that the corresponding  connection forms $\omega$ and $\widetilde{\omega}$ are related by 
$$
  \widetilde{\omega} = \omega - *du.
$$
Differentiating this equation,  using the definition of the Laplacian and the relation   \\  $d\widetilde{A} = \widetilde{\theta}^1 \wedge \widetilde{\theta}^2 = e^{2u}dA$, we obtain 
$$
  d\widetilde{\omega} = d\omega - d*du =  d\omega  + \Delta_g u dA.
$$
\qed

\medskip

We can now prove that any smooth Riemannian metric on a surface is locally conformally equivalent to the flat Euclidean metric: 

\medskip

\begin{definition}
A coordinate system $(x,y)$ on an open domain  $U \subset S$ is  \emph{conformal} for the metric $g$ if there is a function $\rho: U \to \r$ such that 
$$
 g = \rho(x,y) (dx^2+dy^2).
$$
\end{definition}

Such coordinates are also called \emph{isothermal coordinates}.  \index{isothermal coordinates} 
In those coordinates, the Laplacian  is given by 
\begin{equation}\label{FlatLaplacian}
 \Delta_g u = -\frac{1}{\rho(x,y)}\left(\frac{\partial^2 u}{\partial x^2} + \frac{\partial^2 u}{\partial y^2}\right).
\end{equation}

\smallskip 

\begin{proposition}\label{existenceconformcoord}
 On an arbitrary smooth Riemannian surface $(S,g)$, one can introduce  conformal coordinates in  the neighborhood of any  point.
\end{proposition}

\textbf{Proof.}  
From the theory of Elliptic partial differential equations, we know that  any point admits a neighborhood $U$ where we can solve the equation 
 $\Delta_g u = -K$, where  $\Delta_g $ is the Laplacian of $g$ and $K$ its curvature. Applying the previous Lemma, we see that   $\tilde{g} = e^{2u}g$ is a flat metric (that is a metric of vanishing curvature) on $U$ and we can thus find  a system of  Euclidean coordinates $x,y$ in the neighborhood of each point in  $U$ such that  $\tilde{g} = dx^2+dy^2$. In this neighborhood we then have
$$
  g = e^{-2u}(dx^2+dy^2).
$$
\qed

\bigskip

\textbf{Remarks.}  \  \textbf{{(1)}} This result  is classically called the \emph{existence Theorem of isothermal coordinates}. It has an interesting history, going back to Gauss in the case where the metric $g$ is analytic. In 1914 and 1916, Korn and Lichtenchtein  \cite{K,L}  independently proved  this result assuming $g$ to be  only Hölder continuous; and Chern \cite{Chern}  gave a simpler   proof in 1955. In \cite{AhlforsBers}, L. Ahlfors and L. Bers extended this result to the case of surfaces whose metric is only measurable, under certain conditions. 

\smallskip

\textbf{(2)}  Equation (\ref{eq.chgt.curvature})  is usually called the \textit{Liouville Equation}; It plays an important role in differential geometry. 
Observe that integrating this Equation yields:
$$
 \int_S  \tilde{K}d\tilde{A} = \int_S  \tilde{K}e^{2u}dA = \int_S  KdA +  \int_S  \Delta_gu \, dA.
$$
If $S$ is a closed Surface,  we have  $\int_S  \Delta_g udA = -\int_S d * du = 0$ by Stoke's Formula. It follows that 
$$
 \int_S  \tilde{K}d\tilde{A} =   \int_S  KdA,
$$
which is compatible with the Gauss-Bonnet Formula.

\begin{corollaire} \label{cor1.strcomp}
Any smooth Riemannian metric $g$ on an  oriented surface $S$ defines a complex structure on that surface.
\end{corollaire}

\textbf{Proof.}  One can build an oriented atlas on $(S,g)$ with conformal charts. In such an atlas, all changes of
coordinates mappings satisfy the Cauchy-Riemann Equations  and are therefore holomorphic functions.

\qed

\medskip

{\small  One  can also observe that the dual map of the Hodge star defines  an almost complex structure $J : TS \to TS$, 
and this almost complex structure is integrable since $\dim(S) = 2$. The argument proves both the Proposition \ref{existenceconformcoord} and the Corollary \ref{cor1.strcomp}.}

\medskip

We can now quote the Poincaré-Koebe Uniformisation Theorem:
\begin{theorem}\label{th.confconstcurv}
 Any smooth Riemannian metric $g$ on a compact surface is a conformal deformation of a Riemannian  metric $h$ of constant curvature.
\end{theorem}

\textbf{Proof (sketch).}  Let $(S,g)$ be a closed Riemannian surface. Suppose first that $\chi (S) = 0$, then $\int_S KdA = 0$ and we can find a smooth solution $u\in C^{\infty}(S)$ of the equation 
\begin{equation}
 \Delta_g u = -K,
\end{equation}
see e.g. Theorem \ref{th.GreenFunction} below. The equation (\ref{eq.chgt.curvature}) implies then that $h = e^{2u}g$ is a flat metric.

\medskip 

In the case where $\chi (S) < 0$, the argument is similar: solving  the non-linear equation
\begin{equation}\label{eq.Knonlin}
    \Delta_g u - \tilde{K}e^{2u} = 1,
\end{equation}
we obtain the metric $h = e^{2u}g$ of constant curvature  $-1$. Equation (\ref{eq.Knonlin}) can be solved by a variational method 
as was done in 1969 by Melvyn Berger in \cite{B}, see also \cite{Aubin}. 
In the case where $\chi (S) > 0$, the argument fails but the difficulty can be circumvented by introducing a singularity
and solve a linear partial differential  equation instead of the non linear PDE \eqref{eq.Knonlin}, see \cite [th. 4.2]{Troyanov1990}.  
Equation \eqref{eq.Knonlin} can also be solved using  Ricci flow methods, see \cite{CLT,Hamilton}.

\qed  

\medskip

{\small It is sometimes convenient to  normalize the Riemannian metrics   by requiring that $K=-1$ if $\chi(S) < 0$ and 
 $K= +1$ if  $\chi(S) > 0$. If $\chi(S) = 0$ then $K=0$ and we normalize the metric by  requiring $(S,h)$ to have area $1$.
}

\section{Green Kernel and Potential}

As we have seen in the previous section, the curvature measure  of a conformal deformation of a Riemannian surface involves the
Laplacian of the conformal factor. In this section we discuss the inverse of the Laplacian. The following definition will play an 
important role in this context

\medskip

\begin{definition}
Let $u$ be a locally  integrable function defined on a Riemannian surface $(S,h)$  and $d\mu$ be a signed Radon measure on $S$.
We say that $\mu$ is the Laplacian of $u$ \textit{in the weak sense} (or \textit{in the sense of distributions}), and we write $\Delta_h u =_w \mu$,  if for 
any smooth function $\varphi \in C_0^{\infty}(S)$ with compact support, the following equality holds
\begin{equation}\label{eq.wlaplace}
 \int_S  \Delta_h \varphi (x)  u(x) \; dA_h(x) = \int_S  \varphi (y)  d\mu (y).
\end{equation}
\end{definition}

Using integration by part, we immediately observe that for a $C^2$ function  $u$,  the previous identity holds with $d\mu = \Delta_h u dA_h$.

\medskip

Considering first the standard Euclidean metric $h = dx^2 + dy^2$ in a plane domain $U \subset \r^2$, we state the following basic fact:
\begin{lemme}
If the locally integrable function $u$ on the  domain $U \in \mathbb{R}^2$ admits a weak Laplacian $d\mu$ with respect to the Euclidean metric  $dx^2+dy^2$, then $u$ can be written as 
$$
 u(x) = \psi(x) - \frac{1}{2\pi}  \int_U \log |x-y| d\mu (y),
$$
where $\psi$ is a harmonic function.
\end{lemme}

\textbf{Proof.}
Recall that for any smooth function $\varphi$  with compact support on a domain $U \subset \mathbb{R}^2$ we have for any
$y \in U$
$$
  \varphi (y) = -\frac{1}{2\pi}  \int_U \Delta \varphi (x) \log |x-y| dx,
$$
where $\Delta = -\left(\frac{\partial^2}{\partial x^2} +\frac{\partial^2}{\partial x^2}\right)$ is the standard  Laplacian. Therefore if we define
\begin{equation}\label{logpotential}
 v(x) = \frac{1}{2\pi}  \int_U \log |x-y| d\mu (y),
\end{equation}
then $d\mu$ is the weak Laplacian of $v$. Indeed we have
\begin{align*}
 \int_U  v(x)\Delta \varphi (x) dx &=\int_U  \left(  - \frac{1}{2\pi}   \int_U   \log |x-y| d\mu (y) \right) \Delta \varphi (x)  dx
  \\ &=  \int_U   \left( -\frac{1}{2\pi}  \int_U   \log |z-y| \Delta \varphi (x)  dx\right)  d\mu (y)
  \\ &=    \int_U \varphi (y)d \mu (y).
\end{align*}
The Lemma is proved since $\psi = u - v$ is clearly harmonic on $U$.

\qed


\medskip

The function $v$ defined by (\ref{logpotential}) is called the \textit{logarithmic potential of the measure
$d\mu$}. \index{logarithmic potential} An equivalent construction holds on compact surfaces and is based on the following Theorem:

\begin{theorem}\label{th.GreenFunction}
Let $(S,h)$ be a closed smooth Riemannian surface. There exists a unique function 
$G : S\times S \to  \r \cup \{ +\infty\}$ satisfying the following conditions:
\begin{enumerate}[(a)]
  \item $G$ is  $C^{\infty}$ on  $ S\times S \setminus \{(x,x) \tq x \in S \}$;
  \item $G(x,y) = G(y,x)$;
  \item $|G(x,y)| \leq C \cdot (1+ |\log d(x,y)|)$ for some constant $C$;
  \item $\ds \int_S G(x,y) dA_h (y) = 0$;
  \item The following identity holds for any function $u \in C^2(S)$:
\begin{equation}\label{eq.uFrmG}
   u(x) =  \int_S G(x,y) \Delta u (y)  dA_h (y) + \frac{1}{\aire (S)} \int_S  u (y)  dA_h,
\end{equation}
  where  $\Delta = \Delta_h$ is the Laplace-Beltrami operator with respect to the metric $h$.
\end{enumerate}
\end{theorem}
A proof can be  found in  \cite{Aubin} or  \cite{deRham}; the function $G(x,y)$ is called the \textit{Green Kernel} 
\index{Green kernel} of the Riemannian manifold $(S,h)$.

\begin{proposition}\label{prop.GsolvesLapl}
Let $\mu$ be a signed Radon measure on the closed  Riemannian surface  $(S,h)$. If  $\mu$ has vanishing integral on $S$, then the function 
\begin{equation}\label{eq.pot}
  u(x) =  \int_S G(x,y) d\mu (y)
\end{equation}
is a weak solution of $\Delta_h u =_w \mu$.
\end{proposition}

\textbf{Proof.}  We need to check that (\ref{eq.wlaplace}) holds for the function $u$ defined by (\ref{eq.pot}). Indeed we have
\begin{align*}
 \int_S  \Delta \varphi (x) \cdot u(x) \; dA(x) & = \int_S \int_S  \Delta \varphi (x)  G(x,y)  \; d\mu (y)dA(x) 
  \\  & = \int_S \int_S  G(x,y) \Delta \varphi (x)   \; dA(x)  d\mu (y)
  \\  & =  \int_S \left( \varphi (y)  - \overline{\varphi }\right)  d\mu (y)
  \\  & =  \int_S  \varphi (y)  d\mu (y),
\end{align*}
where we have noted 
$\overline{\varphi } =  \frac{1}{\aire (S)} \int_S   \varphi (y)   dA_h (y)$ the average value of $\varphi$. The last equality follows from the hypothesis  $\int_S d\mu = 0$.

\qed

\medskip

\begin{definition}
 The function  $u$  defined by (\ref{eq.pot}) is the \emph{potential} \index{potential (of a measure)}  of the measure $d\mu$ relative to the metric  $h$. 
\end{definition}

\medskip

In what follows, we will denote by $\mathcal{V}(S,h)$ the space of functions $u$ on  $S$  whose weak Laplacian
is a well defined signed Radon measure. Some properties of this space are: 

\medskip

\begin{proposition}\label{propVsh}
For a given closed Riemannian surface $(S,h)$ the following holds:
\begin{enumerate}[(i)]
  \item A function $u$ on $S$ belongs to $\mathcal{V}(S,h)$  if and only if  $u$ is the difference of two subharmonic functions.
  \item The weak Laplacian $\mu$ of a function  $u$ in  $\mathcal{V}(S,h)$ has vanishing integral: $\mu(S) = \int_S d\mu = 0$.
  \item A function  $u$ in  $\mathcal{V}(S,h)$ can be recovered from its weak Laplacian $\mu$ from the following  integral formula: 
\begin{equation}\label{eq.pot2}
  u(x) =  \int_S G(x,y) d\mu (y) +  \frac{1}{\aire (S)} \int_S  u (y)  dA_h.
\end{equation}
  \item For any $1 \leq p < 2$  we have  $\mathcal{V}(S,h) \subset W^{1,p}(S)$  where $W^{1,p}(S)$ is the first order $p$-Sobolev space,  and we have  
$$
 \sup_{0< \varepsilon \leq 1} \sqrt{\epsilon} \|  \nabla u \|_{L^{2-\varepsilon}(S)} \leq C(h, |\mu| (S)).
$$
  \item Any function $u$ in  $\mathcal{V}(S,h)$  is almost everywhere approximately differentiable and we have  
$$
  \| \text{ap} \nabla u \|^*_{L^{2,\infty}(S)} \leq C(h, |\mu| (S)),
$$
where  $L^{2,\infty}(S)$ is the  Lorentz space.
\end{enumerate}
\end{proposition}

\textbf{Proof.} Property (i) follows from the decomposition of the measure $\mu = \Delta_h u$ as a difference of two non negative measures.

Property (ii) follows from applying (\ref{eq.wlaplace}) to the constant function $\varphi = 1$. 

Property (iii)  is contained in Theorem \ref{th.GreenFunction} for the special case of  smooth functions. To prove it in  the general case, consider an arbitrary  function $u \in \mathcal{V}(S,h)$ and  define a function $v$ on $S$ by
$$
  v(x) =  \int_S G(x,y) d\mu (y) +  \frac{1}{\aire (S)} \int_S  u (y)  dA_h
$$
It then follows from Proposition \ref{prop.GsolvesLapl} that  $\mu$ is also the weak Laplacian of $v$. Therefore the weak Laplacian of
$(u-v)$ vanishes and by Elliptic regularity it is a smooth harmonic function. Applying the maximum principle, we see that $(u-v)$
is constant since  $S$ is a closed surface. From  statement (d) in Theorem \ref{th.GreenFunction} we know  that $\int_S (u-v) dA = 0$, this proves the equality $u = v$.

\medskip 

The last two properties will not be used in this paper;  Property (iv) follows from  \cite[Theorem 9.1]{Stampacchia} and \cite[Theorem 2]{Iw} and a proof of Property (v)  can be found in  \cite[Theorem 2]{DHM}.

\qed

\medskip

The following Lemma says that  the weak Laplacian  of a function on a Riemannian surface only depends on the conformal class of the Riemannian metric:

\begin{lemme}  Let $h_1$ and $h_2$ be two conformal metrics on a closed surface $S$. If  $u\in \mathcal{V}(S,h_1)$ 
then $u\in \mathcal{V}(S,h_2)$ and the weak Laplacian of $u$ with respect to $h_2$, seen as a signed measure on $S$, is equal to the 
weak Laplacian of $u$ with respect to $h_1$.
\end{lemme}

\medskip

\textbf{Proof.}  Let us denote  by $\mu_i$  the  weak Laplacian of $u$  with respect to $h_i$ and  by $dA_i$ the area measure of 
$h_i$ for $i=1,2$. Because the Hodge $*$ operator on $1$ forms on a surface is obviously a conformal invariant, we have for 
any smooth function $\varphi$ on $S$:
$$
 \Delta_{h_1} (\varphi) dA_1 = - d*d\varphi   = \Delta_{h_2} (\varphi) dA_2.
$$
Therefore 
$$
 \int_S   \varphi d\mu_1  = \int_S u \ \Delta_{h_1} \varphi dA_1 =  \int_S u \ \Delta_{h_2} \varphi dA_2 
= \int_S   \varphi d\mu_2,
$$
and it follows that $\mu_1 = \mu_2$.

\qed

\medskip

\bigskip

For more on subharmonic functions and potential theory, we refer to the books \cite{AG} and \cite{Tsuji}.

\section{The Conformal Representation of  Alexandrov Surfaces}

We have defined the  Alexandrov metrics with bounded integral curvature on a surface $S$ as limits of sequences of  Riemannian metrics with a uniform 
bound on the total absolute curvature. In this section we present a different viewpoint,  where Alexandrov surfaces are 
seen as conformal deformations of smooth Riemannian surfaces. This point of view is based on the work of  Yuri Reshetnyak and Alfred Huber in the 1950-60's. The basic construction is given in the following construction:
 Let $(S,h)$ be a connected closed Riemannian surface and $u\in \mathcal{V}(S,h)$. For any $x,y \in S$, we define 
\begin{equation} \label{def.qdist}
 d_{h,u}(x,y) = \inf  \left\{ \int_0^1 e^{u(\alpha(t))} \|\dot \alpha(t)\|_hdt  \mid  \alpha \in \mathcal{C}_S(x,y) \right\},
\end{equation}
where $\mathcal{C}_S(x,y)$ is the set of  rectifiable  curves  $\alpha : [0,1] \to S$ such that $\alpha(0) = x$ and $\alpha (1) = y$.  
Note that the set of points where a given function  $u\in \mathcal{V}(S,h)$ is undefined is a polar set (see e.g. \cite[chap. 5]{AG}), in particular it has 
Hausdorff dimension $0$. The integral in (\ref{def.qdist}) is then well defined and it  clear that $d_{h,u}$ is a  pseudometric, that is  
$0 \leq d_{h,u}(x,y) \leq \infty$ for any  $x,y \in S$, and 
$d_{h,u}$ is symmetric and satisfies the triangle inequality.


\medskip 

\begin{theorem}[Reshetnyak] \label{th.msubh}
The pseudometric  $d_{h,u}$ satisfies the following properties:
\begin{enumerate}[(i)]
\item  $d_{h,u}$ is separating, that is $d_{h,u}(x,y) > 0$ if $x\neq y$.
\item  $d_{h,u}(x,y) < \infty$ for every pair of points $x,y\in S$ such that $\mu (\{ x\}) < 2\pi$ and $\mu(\{ y\}) < 2\pi$, where $\mu$ is the measure $\Delta_h u$.
\item A point $x\in S$ such that $\mu (\{ x\}) > 2\pi$ is a point at infinity, that is any non constant curve containing that point has infinite length.
\item A point $x$ such that $\mu (\{ x\}) = 2\pi$ may be at finite or infinite distance from ordinary points.
\end{enumerate}
\end{theorem}

The first assertion  follows from  \cite[Theorem 4.1]{R60b}  and the second one is proved in \cite[Theorem 5.1]{R60b} (see also Lemma 4.1 and 4.2 in  \cite{Troyanov1991}). Statement (iii) follows from Theorems 12.1 and 12.2 in \cite{R60a} and the last statement is seen by considering suitable examples such as Example 5 below.

\qed 

\medskip

\begin{definition}
Following Reshetnyak, we call a metric $d$ on the surface $S$ a \textit{subharmonic metric}  \index{subharmonic metric}  if there exists a Riemannian metric $h$ and 
$u\in \mathcal{V}(S,h)$ such that $d = d_{h,u}$, defined in \eqref{def.qdist}.
\end{definition}
 The terminology refers to the fact that $u$ is a difference of two subharmonic functions.
The next result states that subharmonic metrics on a surface  are precisely the Alexandrov metrics.
\begin{theorem}[Reshetnyak--Huber] \label{th.RH}
 Let $S$ be a connected closed smooth surface. Then the following holds:
\begin{enumerate}[(a)]
  \item If $h$ is a Riemannian metric on $S$ and $u\in \mathcal{V}(S,h)$ are such that the pseudo-metric $d = d_{h,u}$ has no point at infinite distance, then $(S,d)$ is an  Alexandrov surface with bounded integral curvature. 
  \item Conversely,  for any closed surface with  bounded integral curvature  $(S,d)$  in the sense of Alexandrov, there exists a smooth Riemannian metric $h$ on $S$ and a function $u\in \mathcal{V}(S,h)$ such that  $d = d_{h,u}$. \\
\end{enumerate}
\end{theorem}

A local version of this remarkable Theorem has been proved by Reshetnyak in \cite[Theorem I and II]{R60b} and the global version is due to Huber \cite{Huber}.  Reshetnyak's proof is based on a  fundamental convergence theorem for subharmonic metrics on the plane. To state this result, consider  two  sequences of non-negative measures  $(\omega'_n)$ and  $(\omega''_n)$ whose support are contained in 
a fixed compact subset of $\mathbb{C}$. Suppose that these measures weakly converge to $\omega'' = \lim_{n\to \infty} (\omega''_n)$ and
$\omega' = \lim_{n\to \infty} (\omega'_n)$, and define the  following (singular) conformal Riemannian metrics on the plane:
$$
 g_n = e^{2u_n}|dz|^2  \quad \text{and} \quad  g = e^{2u}|dz|^2,
$$
where $u_n$ and $u$ are the logarithmic potential of $\omega_n = \omega'_n -\omega''_n$ and $\omega= \omega' - \omega''$, that is 
$$
  u_n(z) = -\frac{1}{\pi} \int_{\mathbb{C}} \log|z-\zeta| d\omega_n(\zeta)  \quad \text{and} \quad 
   u(z) = -\frac{1}{\pi} \int_{\mathbb{C}} \log|z-\zeta| d\omega(\zeta).
$$
With these notations, Reshetnyak proved the following 
\begin{theorem}[Theorem III in \cite{R60a}] \label{th.III}
If  $\Omega \subset \mathbb{C}$ is a relatively compact domain with piecewise smooth boundary containing no point $z$ such that  $\omega_1(\{z\}) \geq 2\pi$,
then the metric $d_n$ induced by $g_n$ on  $\Omega$ converges uniformly to  the subharmonic  metric $d$ induced by $g$, that is 
$$
  D_{\Omega}(d_n,d)  \to 0.
$$
\end{theorem}

The condition  $\omega_1(\{z\})  < 2\pi$ is here to  keep us safe from the appearance of cusps, in particular  infinite cusps, in the construction.
We now describe the curvature measure of a subharmonic metric. 
\begin{proposition}
The curvature measure of the Alexandrov surface $(S,d_{h,u})$ is given by
\begin{equation}\label{eq.chgt.courbure2}
 d\omega = K_hdA_h + d\mu,
\end{equation}
where $K_h$ is the Gauss curvature of $h$, $dA_h$ is the area measure and $\mu =_w \Delta_h u$ 
is the weak Laplacian of $u$ for the metric $h$.
\end{proposition}

\textbf{Proof.}  If $u$ is smooth, then the Proposition has been established earlier in Lemma \ref{lemconformalcurvature}.
For the general case, one can introduce local isothermal coordinates in the neighborhood of any point in $S$ and apply Theorem \ref{th.III}.

\qed

\medskip

\textbf{Remark.}  Observe that the Gauss-Bonnet  formula (\ref{eq.GB}) can be recovered from (\ref{eq.chgt.courbure2}) and statement (ii) in  Proposition  \ref{propVsh}, since
$$
\int_S d\omega  = \int_S K_hdA_h + \int_Sd\mu   = \int_S K_hdA_h= 2\pi \chi (S).
$$

 \medskip
 
The next result states that a subharmonic metric is associated to a unique  conformal structure on the surface.
 
  \medskip

\begin{theorem}\label{th,rigiditeconforme}
Any isometry between two surfaces with subharmonic metrics is a conformal map. More precisely,   
if 	$(S,h)$  and  $(S',h')$ are two closed Riemannian surfaces  and   $d= d_{h,u}$ and $d' =d_{h',u'}$ are conformal  subharmonic 
metrics, then any  distance preserving map  $f : (S,d) \to (S',d')$  is a  conformal map from $(S,h)$ to  $(S',h')$.
\end{theorem}

This result is Theorem 7 in \cite{R63}. Reshetnyak's  proof is based on a 1937 Theorem by D. E. Menchoff which says   that a $1$-quasiconformal map between two  Riemann surfaces is a conformal map (see \cite{G,Golberg,M}). Huber gave a different proof  in \cite[Lemma 7]{Huber}.

\medskip

\begin{corollaire}
A metric $d$ with  bounded integral curvature on an oriented surface $S$  is associated to a unique complex structure on that  surface.
\end{corollaire}

 \textbf{Proof.}  
 The previous Theorem tells us that $d = d_{h,u}$ where $h$ is a Riemannian metric on $S$ and $u\in \mathcal{V}(S,h)$. By Corollary 
\ref{cor1.strcomp} we know that $h$ defines a complex structure on $S$ and from 
Theorem \ref{th,rigiditeconforme}  we know that this complex structure   is uniquely determined from the metric $d$.

\qed

\medskip 

\textbf{Examples.} Let us now give some concrete examples of surfaces with bounded integral curvature in their conformal representation: 

\smallskip 

(1) \  Let $V$ be an Euclidean cone of total angle $\theta$. Then $V$ is isometric to  $\mathbb{C}$ with the metric
$$
  ds^2 = |z|^{2\beta}|dz|^2 \, ,
$$
where  $\beta = (\frac{\theta}{2\pi} -1)$ (see \cite[prop.1]{Troyanov1986}).This metric is of class  $L^p_{loc}$ for any  $1 < p < -1/\beta$ \  if $\beta < 0$ (that is if $\theta < 2\pi$) and of class  $L^{\infty}_{loc}$ if $\beta > 0$. Its curvature is the measure 
$$d\omega = -2\pi \beta \cdot \delta _0,$$
where $\delta _0 $  is the Dirac measure  at $0$.
\medskip 

(2) \ Consider the surface $S_0$ obtained by gluing a half-sphere of radius $1$ to a half-cylinder of same radius along their boundaries
by an isometry.  Then $S_0$ is isometric to   $(\mathbb{C},ds^2)$ where $ds^2 = \rho (z) |dz|^2$ is given by 
$$
\rho (z) = \begin{cases}
\tfrac{4}{(1+|z|^2)^2}  & \text{if  $|z| \le 1$,} \\ \ds
\tfrac{1}{|z|^2}  & \text{if  $|z| \ge 1 $.} 
\end{cases}
$$
This metric is of class $C^{1,1}$ and its curvature measure is absolutely continuous, given by 
$d\omega = KdA$,  where 
$$
K (z) = \begin{cases}
1 \,  & \text{if  $|z| < 1$,}\\
0 \,   & \text{if  $|z| > 1$.} 
\end{cases} 
$$
One may consider the  one-point compactification of $S_0$, let us denoted it by $S = S \cup \{\infty\}$. Then $S$ is a topological sphere, conformally equivalent to $\mathbb{C} \cup \{\infty\}$, and the curvature measure is now given by 
$$
 d\omega = KdA + 2\pi \delta_{\infty}.
$$

\medskip

(3)   \ Let $S$  be the surface obtained by gluing two Euclidean discs of radius $1$ along their boundaries by an isometry. 
Then  $S$ is a topological sphere with a singular metric which  is flat on the complement of a singular circle $\Sigma$. 
This surface is isometric to   $\mathbb{C} \cup \{ \infty\}$ with the metric $ds^2 = \rho (z) |dz|^2$ where 
$$
\rho (z) =  \begin{cases}
1 \,  & \text{if  $|z| \le 1$,} \\ \ds
|z|^{-4} \,   & \text{if  $|z| \ge 1$.} 
\end{cases} 
 $$
This metric is of class $C^{0,1}$ (Lipschitz) and its curvature measure is 
$$d\omega = 2 \cdot ds _{| \Sigma},$$
where $ds _{| \Sigma}$  is the length measure along $\Sigma$.

\medskip

(5) \ (From Hulin-Troyanov \cite{HT}). The metric in  $\{z \in \mathbb{C} \mid |z| < 1\}$ defined by
$$
 g = \frac{|dz|^2}{|z|^2 |\log |z||^{2a}}
$$
has  Gauss curvature  $K = - a |\log |z||^{2a-2}$ on $\{z\neq 0\}$ and the curvature at the origin is  $\omega(\{0\}) = 2\pi$.
The curvature measure is then 
$$
  d\omega = 2\pi \delta_{0} - a |\log(|z|)|^{2a-2}dA.
$$
If $a > 0$, the singularity at the origin is a cusp and is $a > 1$, the cusp is at finite distance from a regular point 
(say the point $z = 1/2$). If  $0 < a \leq 1$,
then the  cusp is at infinite distance (any non trivial curve containing the origin has infinite length). 
The special case $a=1$ is  the  Beltrami pseudosphere and the limit case $a = 0$ is a half-cylinder.

\section{Toward a classification of compact Alexandrov Surfaces}
\label{sec.class}

So far we have associated the following data to any  Alexandrov metric without cusp on a closed surface $S$:
\begin{enumerate}[(i)]
 \item A well defined conformal structure on $S$.
 \item A signed Radon measure $d\omega$ such that $\omega (\{x\}) < 2\pi$ for any point $x$ on $S$ and  the Gauss-Bonnet 
 condition $\int_S d\omega = 2\pi \chi (S)$ holds.
\end{enumerate}

In the converse direction one may state the following Theorem, which basically summarizes the previous results: 
\begin{theorem}\label{th.clsf}
For any conformal structure on a closed surface $S$  and any signed Radon measure  $d\omega$  on  $S$ such that 
$\int_S d\omega = 2\pi \chi (S)$ and  $\omega (\{x\}) < 2\pi$  for any point  $x\in S$,  there exists an Alexandrov metric on $S$ 
belonging to that conformal structure and whose curvature measure is  $d\omega$.
This metric is unique up to a possible dilation.
\end{theorem}

An explicit example of surface with fractal curvature is discussed in the recent paper \cite{Dima}.

\medskip 

\textbf{Proof.}  To prove the existence of a compatible Alexandrov metric, we first fix a Riemannian metric $h$ with constant curvature $K_h$ on $S$ in the given conformal class.  Set $d\mu =  d\omega - K_hdA_h$ and observe that from Gauss-Bonnet formula, we know that 
$d\mu$ has vanishing integral on $S$.  Choosing  the function   $u\in \mathcal{V}(S,h)$ to  be the potential of the measure $d\mu$,  the desired metric is given by $d = d_{h,u}$. 
 
\smallskip

To prove uniqueness, we consider another Alexandrov metric  $d'$ on $S$ in the same conformal class and with the same curvature measure. By  Theorem  \ref{th.RH},  there exists  a smooth metric $h'$ on  $S$  and a function $u'\in \mathcal{V}(S',h')$  such that $d' = d_{h',u'}$. Theorem  \ref{th,rigiditeconforme} implies that $h$ and $h'$ are conformally equivalent;  there exists therefore a smooth function  $v \in C^{\infty}(S)$ such that $h'  = e^{2v}h$.  Replacing $u'$ by $u'+v$ if necessary, we may assume that  $h=h'$. We therefore have  $d' = d_{h,u'}$  and its curvature measure is $d\omega$. Then 
$$
 \Delta_h u' = d\omega  - K_hdA_h =  \Delta_h u,
$$
that is $\Delta_h (u'- u) = 0$ and thus  $(u'-u)$  is constant. This implies that the metric  $d'$ is a (constant) multiple of $d$.

\qed 

\medskip

A first consequence is the following result classifying  Euclidean surfaces with conical singularities (see \cite{Troyanov1986}).

\begin{corollaire}
Let  $S$ be a closed surface $x_1,\cdots, x_n$ be points on  $S$ and $\theta_1, \cdots ,\theta_n > 0$. Suppose that  $\sum_i(2\pi-\theta_i)=2\pi\chi(S)$, then for any conformal structure on $S$, there exists a conformal polyhedral metric on $S$ with a conical singularity of angle  $\theta_i$ at $x_i$ for  $(i=1, \cdots,n)$. This metric is unique to homothety.
\end{corollaire}

\textbf{Proof.}  Apply the previous Theorem to the discrete measure $d\omega = \sum_i(2\pi-\theta_i) \delta_{x_i}$, where $\delta_{x}$ is the Dirac mass at $x \in S$.

\qed

\bigskip

\textbf{Remark.} Theorem \ref{th.clsf} can  be seen as a classification Theorem for Alexandrov surfaces. 
Let $S$ be an oriented closed surface and let $\mathcal{M}_0(S)$ denote the space of metrics with bounded integral curvature on $S$ 
without cusps, $\mathcal{C}(S)$ be the space of conformal structures
and $\mathcal{R}_{2\pi}(S)$  be the space of signed  Radon measures  $d\omega$  on $S$ such that
$$
 \int_S d\omega = 2\pi\chi(S) \quad \text{and} \quad \omega(\{x\}) < 2\pi
 \ \text{ for all } x \in S.
$$
Then the previous Theorem states that we have a natural bijection 
\begin{equation}\label{eq.clsf}
 \mathcal{M}_0(S)  \cong   \mathcal{C}(S) \times \mathcal{R}_{2\pi}(S)\times \r_+
\end{equation}
(the factor $\r_+$ controls the scaling of the metric). 

\bigskip

Combining   Theorem \ref{th.clsf}   with Alexandrov's solution of  Weyl's problem (see Section \ref{WeylPb}), we obtain
the following 

\begin{corollaire}\label{covexdomain}
Any non negative Radon measure $d\omega$ on  $S^2$ such that  $\int_{S^2}d\omega = 4\pi$ and  $\omega (\{x\}) < 2\pi$ for all $x \in S^2$  is the intrinsic curvature measure on the boundary of a bounded convex domain $\Omega \subset \r^3$.
Two such measures  $d\omega_1$ and $d\omega_2$ are realized by two similar  convex domains $\Omega_1$ and $\Omega_2$  if and only if
there is a conformal transformation of the sphere transforming $d\omega_1$ to $d\omega_2$.
\end{corollaire}

Recall that two domains  $\Omega_1$ and $\Omega_2$ in $\r^3$ are \textit{similar} if there is an isometry $f : \r^3 \to \r^3$ and a positive  constant $\lambda$ such that $\lambda f(\Omega_1) = \Omega_2$.

\section{Some Questions and Problems}

We conclude   with a short list of questions and open problems related to the results in this paper:
 
\bigskip 
 
\textbf{Problem 9.1.}
Can one prove the following  global version of  Reshetnyak's convergence Theorem:   \ 
\emph{Let  $(S,h)$ be a smooth Riemannian surface and  
$(d\mu^+_n)$, $(d\mu^-_n)$ 
be two sequences of  non negative Radon measures on $S$ that  weakly converge respectively to the  measures $d\mu^+ = \lim_{n\to \infty} d\mu^+_n$ and  $d\mu^- = \lim_{n\to \infty} d\mu^-_n$. 
Suppose $\mu^+ (\{ x\}) < 2\pi$ for any  point $x\in S$ and let us denote by  $u_n$ the potential  of $d\mu_n = d\mu^+_n - d\mu^-_n$ and by  $u$  the potential of $d\mu = d\mu^+ - d\mu^-$. Can we then conclude that  
$d_{h,u_n} \to d_{h,u}$ in the uniform topology}\footnote{
This statement  was somewhat boldly stated as  Theorem 6.2 in \cite{Troyanov2009}. When I wrote that paper, I  presumed  that following Reshetnyak's arguments should easily lead to the proof. But this is in fact not so obvious and it seems more appropriate to consider it an open problem.
}?

\bigskip 

In relation with this problem, one can cite the work of Clément Débin, who proved in his thesis a global compactness Theorem under some additional geometric bounds (but without fixing the conformal class), see \cite{Debin}.  Note also that limit spaces of Alexandrov's surfaces with uniformly bounded diameter
and total absolute  curvature on a closed surface under Gromov-Hausdorff  convergence have been described in \cite{Shioya}. These limits spaces are not necessarily topological surfaces.

\bigskip 

\textbf{Problem 9.2.} \textit{Are there some natural topologies on the spaces $\mathcal{R}_{2\pi}(S)$  and  $\mathcal{C}(S)$ for which  the bijection  (\ref{eq.clsf})  is  a homeomorphism  when $\mathcal{M}_0(S)$ is equipped with the uniform distance?}

\medskip

{Recall that the $*$-weak topology on the space of Radon measure is  not metrizable, but a natural metric compatible with the weak convergence
can be defined on  $\mathcal{R}_{2\pi}(S)$, see e.g. \cite[page 132]{Doob}.
}

\medskip

{\small Note also that the  group of homeomorphisms  of the surface acts  naturally on both sides of  the bijection  (\ref{eq.clsf}). \textit{We then also ask as a subproblem  if   the map   (\ref{eq.clsf})  defines a homeomorphisms of the corresponding quotients spaces.}  Recall that  the quotient space of $\mathcal{C}(S)$ by the  group of  homeomorphisms of $S$  is the familiar \textit{moduli space} of the surface $S$, well understood from Teichm\"uller theory.}

\bigskip 

\textbf{Problem 9.3.}  \textit{Try and extend the theory to include the cusps (at least the finite cusps), which are natural metrics sitting on the boundary of $\mathcal{M}_0$.}

\bigskip 

\textbf{Problem 9.4.} Theorem \ref{covexdomain} states that there is a map from the set of non negative measures on $S^2 = \mathbb{C} \cup \{\infty\}$ such that $\omega(S^2)
 = 4\pi$  (up to Möbius tranformations) to the set of bounded convex domains in $\r^3$ (up to similarity), the map being such that the 
given measure corresponds to the curvature of the boundary of the convex domain. \textit{We may then ask to describe this map as explicitly as possible and study its properties.}

\bigskip 

\textbf{Problem 9.5.} 
Recall that the \textit{Hilbert metric}   \index{Hilbert metric} $h$ is the distance defined on a given bounded convex domain $\Omega \subset \r^3$  as 
$$
  h(x,y) = \frac{1}{2} \log \left(\frac{|a-x| |y-b|}{|b-x| |a-y|} \right),
$$
for any $x,y \in \Omega$, where $a$ and $b$  are   the intersections of the line through $x$ and $y$ with $\partial \Omega$  and $y$ lies between $x$ and $a$. The Hilbert  distance is invariant under projective transformations leaving the domain  $\Omega$ invariant, it is a complete Finslerian  metric. 
\textit{We then ask to investigate the possible relations between the Hilbert geometry of $\Omega$ and the curvature measure of $S = \partial\Omega$.} As a subproblem one could consider the case where the curvature measure $d\omega$ on $S$ is singular (i.e.  has no absolutely continuous part) and try to relate the Hausdorff dimension of the support of $d\omega$  to the volume entropy of 
$(\Omega, h)$. The paper \cite{Vernicos2017} might provide some useful hints. See also  \cite{Merlin} for a related discussion in dimension $2$.

\bigskip 

\paragraph{Acknowledgement.}
The author is thankful to François Fillastre for the invitation to translate and update the original paper, for carefully reading the manuscript and
his precise and useful comments. He also thanks Boris Buffoni for fruitful discussions related to some analytic aspects in this work. 
This research was in part supported by Swiss SNF grant  200021L--175985.
 

 \printindex 

\end{document}